\documentclass[onefignum,onetabnum]{siamonline220329}

\usepackage{lipsum}
\usepackage{amsfonts}
\usepackage{graphicx}
\usepackage{epstopdf}
\usepackage{algorithmic}
\ifpdf
\DeclareGraphicsExtensions{.eps,.pdf,.png,.jpg}
\else
\DeclareGraphicsExtensions{.eps}
\fi

\usepackage{enumitem}
\setlist[enumerate]{leftmargin=.5in}
\setlist[itemize]{leftmargin=.5in}

\newsiamremark{remark}{Remark}
\newsiamremark{hypothesis}{Hypothesis}
\crefname{hypothesis}{Hypothesis}{Hypotheses}
\newsiamthm{claim}{Claim}
\newsiamremark{example}{Example}

\headers{A discrete Kolmogorov  predator-prey model}{Lei Niu and Susu Wang}

\title{Global stability and period-doubling bifurcations of a discrete Kolmogorov predator-prey model with Ricker-type prey growth\thanks{The work was supported by the National Natural Science Foundation of China (No. 12001096) and the Fundamental Research Funds for the Central Universities (No. 2232023A-02 and 2232024G-13).}}

\author{
	Lei Niu\thanks{School of Mathematics and Statistics,\! Donghua University,\! Shanghai 201620,\! China (\email{lei.niu@dhu.edu.cn}).}
	\and
	Susu Wang\thanks{School of Mathematics and Statistics,\! Donghua University,\! Shanghai 201620,\! China (\email{wangsusu@mail.dhu.edu.cn}).}
}

\usepackage{amsopn}

\ifpdf
\hypersetup{
  pdftitle={Global stability and period-doubling bifurcations of a discrete Kolmogorov predator-prey model with Ricker-type prey growth},
  pdfauthor={Lei Niu and Susu Wang}
}
\fi

\usepackage{amsmath}
\allowdisplaybreaks[4]
\usepackage{float}
\usepackage{subcaption}

\begin{document}

\maketitle
\begin{abstract}
In this paper, we study the dynamics of a  discrete Kolmogorov predator-prey model with Ricker-type prey growth.  
We give the sufficient and necessary condition to guarantee the existence and uniqueness of the positive fixed point. Using the center manifold theory, we prove that the period-doubling bifurcations can occur at the positive fixed point. Furthermore, our numerical simulations reveal that the model can exhibit cascades of period-doubling bifurcations leading to chaos, which is a significant difference from the behavior of continuous predator-prey models. Despite the complexities of the model dynamics, we are able to provide a criterion for the global stability of the positive fixed point by using a geometric analysis of the nullclines.
\end{abstract}



\section{Introduction}
The predator-prey interaction within populations is widely recognized as one of the most interesting areas of inquiry in the fields of biology and ecology. Efforts to understand and clarify predator-prey interactions have been advanced through the development and analysis of mathematical models. The history of the predator-prey models can be traced back to the early 20th century when Lotka \cite{Lotka1925ElementsOP} and Volterra \cite{Volterra1926Variazioni} independently developed mathematical models by using differential equations to describe these interactions, which known as the Lotka-Volterra predator-prey models has since become an iconic model of mathematical biology. Gause \cite{Gauze1934TheSF} further proposed another system of much more general equations which refined the Lotka-Volterra models. Kolmogorov \cite{Kolmogorov1936Sulla} provided a most general system of differential equations, often called Kolmogorov system, which was led to postulate some minimal conditions being satisfied in any realistic predator-prey interaction. By utilizing mathematical frameworks developed by Gause and Kolmogorov, researchers have been able to simulate and predict the outcomes of predator-prey interactions under various conditions, contributing to a deeper understanding of ecological dynamics (see, for example, Albrecht et al. \cite{Albrecht1973}, Collings \cite{Collings1997}, Freedman \cite{Freedman1980}, Hirsch and Smale \cite{Hirsch2012},  Hofbauer and Sigmund \cite{hofbauer1988}, Hsu and Huang\cite{Hus1995}, Kuang \cite{Kuang1990}, May \cite{May1972Limit,May1974Stability}, Ruan and Xiao \cite{RUan2001}).

In the last few years numerous models based on difference rather than differential equations for ecological systems have appeared in the literature. When population numbers are small or when births and deaths occur at discrete times, such as within specific intervals like a generation, discrete models offer a more realistic representation of predator-prey dynamics (\cite{Freedman1980,Gyllenberg1997,Hassell1978,Kot2001}). Many other species, such as monocarpic plants and semelparous species (see Kot \cite{Kot2001} for more example), have discrete nonoverlapping generations and their births occur in regular breeding seasons. In fact, the pioneer discrete-time predator-prey model was proposed in 1930s by Nicholson and Bailey \cite{Nicholson1935}, known as the  Nicholson-Bailey model, which is used to describe the interactions between a population of herbivorous arthropods and their insect parasitoids.  Beddinton et al. \cite{Beddington1975} further modified the Nicholson-Bailey model by introducing the inclusion of density-dependent self-regulation by the prey, in which the prey growth was assumed to be described by the Ricker model \cite{Ricker1954Stock} in the absence of predators. Unlike the continuous-time cases, for which the Poincar\'e-Bendixson theorem holds, the two-dimensional discrete-time predator-prey models may have much more complex dynamics. For example, the modified Nicholson-Bailey model studied by Beddinton et al. in \cite{Beddington1975,Beddington1978} can possess period-doubling bifurcations or Neimark-Sacker bifurcations and even chaos. Therefore, it remains a very challenging topic to study the global behavior of the discrete-time predator-prey systems.

Recently, Ackleh et al. \cite {ackleh2019persistence,ackleh2020Long} proposed a Kolmogorov-type discrete-time predator-prey model of much more general equations from first principles, which is given by 
\begin{equation}\label{model_ackleh}
    \begin{cases}
    x_{k+1}=\phi(x_k)(1-f(y_k)y_k)x_k,\\
    y_{k+1}=sy_k+b(x_k)x_kf(y_k)y_k,
    \end{cases}
\end{equation}
where $x$ and $y$ are the densities of the prey and predator, respectively. The constant
$0<s<1$ denotes the survival of predator, while the three positive functions $\phi$, $f$ and $b$ represent density-dependent prey growth, predator consumption of prey, and the conversion of consumed
prey into new predators. In particular, $f$ denotes the probability that an individual prey is consumed by an individual predator so that it is assumed to satisfy $0<f(y),f(y)y<1$. Usually the prey experiences negative density effects from both prey and predator populations, so it is assumed that $\phi(x)$ is decreasing with respect to $x$ and $f(y)y$ is increasing with respect to $y$. Moreover, 
by considering that the biological restrictions on predator fecundity may be imposed, the predator fecundity $b(y) y$ is usually
assumed to be an increasing function of prey density that saturates at a maximal fecundity. Typical functions that satisfy these assumptions are the Beverton-Holt functions (\cite{Beverton1957})
\begin{equation}\label{equ:funsb-f}
    b(x)=\frac{b_0}{1+\gamma x},~~f(y)=\frac{c}{1+cy},
\end{equation}
where $b_0,\gamma>0$ and $0<c<1$. 
We refer the readers to \cite {ackleh2019persistence,ackleh2020Long} for more details on the underlying biological assumptions.

The system \eqref{model_ackleh} has the advantage of keeping all solutions with positive initial conditions always positive. Ackleh et al. \cite {ackleh2019persistence,ackleh2020Long} studied the system \eqref{model_ackleh}  by assuming that the prey grows according to the Beverton-Holt model in the absence of predators, that is, $$\phi(x)=\frac{r}{1+mx}$$ 
with $r,m>0$, and $b$ and $f$ are given by \eqref{equ:funsb-f}. They proved that the positive fixed point is always unique and asymptotically stable when it exists, and provided sufficient conditions to guarantee the global stability of the positive fixed point by using a Lyapunov function. Ackleh et al. \cite{ackleh2021nullcline} further improved the results on the global stability in \cite{ackleh2020Long} by using an interesting method of geometric analysis of the nullclines.

Motivated by Beddinton et al. \cite{Beddington1975,Beddington1978} and May \cite{May1974}, in this paper we will study the dynamical behavior of the system \eqref{model_ackleh} with $b$ and $f$ given by \eqref{equ:funsb-f}  when the prey growth is described by the Ricker model, i.e.,
\begin{equation}\label{equ:fun-phi}
\phi(x)=e^{r-x}.
\end{equation}
Under these assumptions, the system \eqref{model_ackleh} is written as
\begin{equation}\label{pp_model}
    \begin{cases}
        x_{k+1}=x_ke^{r-x_k}(1-\frac{cy_k}{1+cy_k}),\\
        y_{k+1}=sy_k+\frac{b_0x_k}{1+\gamma x_k}\frac{cy_k}{1+cy_k},
    \end{cases}
\end{equation}
where $b_0,\gamma>0$, $0<c<1$ and $0<s<1$. We prove that the system \eqref{pp_model} can have period-doubling bifurcations (Theorem \ref{thm:fild-stable}) and numerically find that cascades of period-doubling bifurcations can occur which lead to chaos eventually. However, unlike the modified Nicholson-Bailey model (\cite{Beddington1975,Beddington1978}), system \eqref{pp_model} does not have Neimark-Sacker bifurcations. 
Moreover, the system \eqref{pp_model} has revealed fascinating dynamics that differ from the system \eqref{model_ackleh} with the prey growth of Beverton-Holt type studied in \cite{ackleh2019persistence,ackleh2020Long,ackleh2021nullcline}, in which period-doubling bifurcations cannot occur, and chaos has not been found so far. 
Although the dynamics of system \eqref{pp_model} can be very complex, we still provide a criterion (Theorem \ref{thm:global-stability}) that the positive fixed point is globally asymptotically stable by using the geometric method of analyzing the nullclines motivated by Ackleh et al. \cite{ackleh2021nullcline}. Our criterion can be seen as a generalization of their result on systems with monotone prey dynamics in systems with non-monotone prey dynamics.

The paper is organized as follows. Section \ref{sec:extinct} includes some preliminaries, where we prove that system \eqref{pp_model} has a compact absorbing set such that every orbit will enter into it. In Section \ref{sec:fix-point}, we study the existence and stability of the fixed points. In particular, we show that the positive fixed point is unique when it exists and provide a criterion on the global stability of the positive fixed point by using a geometric method. In Section \ref{sec:PD}, we study the bifurcations at the positive fixed point. Using the center manifold theory, we provide the rigorous computation of the period-doubling bifurcations and we numerically find that as the bifurcation parameter varies a cascade of period-doubling bifurcations can occur which lead to chaos eventually. The paper ends with a discussion in Section \ref{sec:conclusion}.

\section{Preliminaries}\label{sec:extinct}
Denote the map $G:[0,+\infty)\to [0,+\infty)$ by 
\begin{equation}
    G(x,y)=(G_1(x,y),G_2(x,y))=\big(xg_1(x,y),yg_2(x,y)\big)
\end{equation}
with $(x,y)\in [0,+\infty)$, 
where 
$$g_1(x,y)=e^{r-x}\frac{1}{1+cy}$$
and
$$g_2(x,y)=s+\frac{cb_0x}{(1+\gamma x)(1+cy)}.$$
Then the system \eqref{pp_model} can be written as 
\begin{equation}
    (x_{k+1},y_{k+1})=G(x_k,y_k).
\end{equation}
Let $\mathcal{K}=[0,K_1]\times [0,K_2]$ with $K_1=e^{r-1}$ and $K_2=\frac{b_0}{\gamma(1-s)}+1$.
\begin{lemma}[Absorbing set]\label{lem:absorption}
The set $\mathcal{K}$ is an absorbing set for system \eqref{pp_model}, that is, $G(\mathcal{K})\subset \mathcal{K}$ and for all $(x_0,y_0) \in  [0,+\infty)^2$, there exists $k_0\in \mathbb{N}$ such that $G^k(x_0,y_0)\in \mathcal{K}$ for all $k\geq k_0$.

\end{lemma}
\begin{proof}
Let $\{(x_k, y_k)\}$ be an orbit of the system \eqref{pp_model} with the initial value $(x_0,y_0)\in[0,+\infty)^2$, that is, 
$$(x_k, y_k)=G^k(x_0,y_0),~~ k=1,2,\ldots.$$
By the definitions of $f$ and $\phi$ in \eqref{equ:funsb-f} and \eqref{equ:fun-phi}, one has 
\begin{equation}\label{equ:leqs_Ricker}
    x_{k+1}=\phi(x_k)x_k(1-f(y_k)y_k)\leq\phi(x_k)x_k,
\end{equation}
which implies that $x_{k+1}\leq e^{r-1}$ for all $k\geq 0$. Then one has 
\begin{equation}\label{equ:xk}
    x_k\in[0,K_1],~~\forall k\geq 1.
\end{equation}
Since $f(y)y<1$ and $b(x)x\leq \frac{b_0}{\gamma}$ for all $x\geq 0$ and $y\geq 0$,  
one has
\begin{equation}\label{equ:y-g}
    y_{k+1}=s y_k+b(x_k)x_k f(y_k)y_k\leq s y_k+\frac{b_0}{\gamma}.
\end{equation}
Consider the sequence $\{g_k\}$ satisfying the system of linear difference equation
\begin{equation}
    g_{k+1}=sg_k+\frac{b_0}{\gamma},~k=0,1,2,\ldots
\end{equation}
with $g_0=y_0$. Since $0<s<1$, it is clear that 
$${\lim_{k\to+\infty}}g_k=\frac{b_0}{\gamma(1-s)}.$$
It then follows from \eqref{equ:y-g} that 
$$\limsup_{k\to+\infty} y_k\leq {\lim_{k\to+\infty}}g_k \leq \frac{b_0}{\gamma(1-s)}.$$ 
Therefore, there exists $k_0\geq 1$ such that $y_k\in[0,K_2]$ for all $k\geq k_0$, and hence $G^k(x_0,y_0)\in \mathcal{K}$ by \eqref{equ:xk} for all $k\geq k_0$. On the other hand, for all $(x_0,y_0)\in \mathcal{K}$, one has
$x_1 \leq K_1$ and 
$$
y_1\leq s K_2+\frac{b_0}{\gamma} < K_2,
$$
by the above arguments, which means that  $G(\mathcal{K})\subset \mathcal{K}$. This completes the proof.
\end{proof}

 \begin{lemma}\label{lemma:leqs_Ricker}
 Assume that $r\leq 1$. Then for any $(x_0,y_0)\in \mathcal{K}$, one has that 
$x_k\leq z_k$ for all $k\geq 0$, where $\{z_k\}$ is the sequence of the one dimensional Ricker model
\begin{equation}\label{equ:1D-Ricker}
    z_{k+1}=\phi(z_k)z_k=e^{r-z_k}z_k,~k=0,1,2,\ldots
\end{equation}
with initial condition $z_0=x_0$.
 \end{lemma}
\begin{proof}
If $r\leq 1$, then $K_1\leq 1$, and hence $\phi(x)x$ is an increasing function from $[0,K_1]$ into $[0,K_1]$. 
Note that whenever $x_k\leq z_k$ one has that
\[  x_{k+1}\leq \phi(x_k)x_k \leq \phi(z_k)z_k = z_{k+1}.\]
It then follows that $x_k\leq y_k$ for all $k\geq 0$ since $x_0=z_0$.
\end{proof}

\section{Fixed points}\label{sec:fix-point}
In this section, we study the existence and stability of the fixed points of system \eqref{pp_model}.
\subsection{Boundary fixed points}
\begin{lemma}\label{lem:extinct}
    The trivial fixed point $(0,0)$ of system \eqref{pp_model} is globally attracting if $r\leq 0$ while unstable if  $r>0$.
\end{lemma}
\begin{proof}
 It is clear that $(0,0)$ is unstable by noticing that $e^{r}$ and $s$ are the two eigenvalues of $DG(0,0)$ if $r>0$. Now assume that $r\leq 0$. Note that $z_k\to 0$ as $k\to +\infty$ for $r\leq 0$, where $\{z_k\}$ is the sequence of the one dimensional Ricker model
\eqref{equ:1D-Ricker}
with $z_0=x_0$. Then by Lemma \ref{lemma:leqs_Ricker}, one has that 
$$
\limsup_{k\to+\infty}x_k\leq \lim_{k\to+\infty}z_k=0.
$$
Thus, $x_k \to 0$ as $k\to+\infty$ 
because $x_k\geq 0$.
It then follows that 
$$\lim_{k\to+\infty}b(x_k)x_k f(y_k) = 0,$$
by the boundedness of functions $b$ and $f$. Choose $\varepsilon>0$ sufficiently small such that  $\varepsilon<1-s$. Then there exists $k_0\geq 0$ such that $b(x_k)x_kf(y_k)<\varepsilon$ for all $k\geq k_0$, and hence
\[y_{k+1}\leq (s+\varepsilon)y_k\]
for all $k\geq k_0$. Thus, $y_k \to 0$ as $k\to+\infty$ because  $s+\varepsilon<1$. This completes the proof.
\end{proof}
\begin{lemma}\label{lem:predator-free}
The system \eqref{pp_model} has a predator-free fixed point $(r,0)$ if $r>0$, and moreover, it is locally asymptotically stable if $0<r<2$ and $s+\frac{cb_0r}{1+\gamma r}<1$.
\end{lemma}
\begin{proof}
It is easy to check that $G(r,0)=(r,0)$, so $(r,0)$ is a predator-free fixed point of $G$ if $r>0$.  Moreover, since the Jacobian matrix
\[DG(r,0)=\begin{pmatrix}
		1-r & -cr\\ 0 & s+\frac{cb_0r}{1+\gamma r}
\end{pmatrix},\]
one has that the predator-free fixed point $(r,0)$ is locally asymptotically stable if
\[ 0<r<2 ~ \text{and} ~ s+\frac{cb_0r}{1+\gamma r}<1. \]
\end{proof}

\begin{theorem}\label{thm:predatorfree-globally}
The predator-free fixed point $(r,0)$ of system \eqref{pp_model} is globally asymptotically stable if $0<r\leq 1$ and $s+\frac{cb_0r}{1+\gamma r}<1$. 
\end{theorem}
\begin{proof}
If $0<r\leq 1$ and $s+\frac{cb_0r}{1+\gamma r}<1$, then the predator-free fixed point $(r,0)$ is locally asymptotically stable by Theorem \ref{lem:predator-free}. Now we show that $(r,0)$ is globally attracting. By Lemma \ref{lem:absorption}, it suffices to prove that 
$$
G^k(x_0,y_0)\to (r,0)
$$
as $k\to +\infty$ for any $(x_0,y_0)\in \mathcal{K}$. Note that $z_k\to r$ as $k\to +\infty$ for $0<r\leq 1$, where $\{z_k\}$ is the sequence of the one-dimensional Ricker model
\eqref{equ:1D-Ricker}
with $z_0=x_0$ (see \cite{May1974}). Then by Lemma \ref{lemma:leqs_Ricker}, one has that 
\begin{equation}\label{equ:upper_bound}
    \limsup_{k\to+\infty}x_k\leq \lim_{k\to+\infty}z_k=r.
\end{equation}
Since $s+\frac{cb_0r}{1+\gamma r}<1$, there exists $\varepsilon>0$ such that
\[
s+\frac{cb_0(r+\varepsilon)}{1+\gamma(r+\varepsilon)}<1.
\]
Then by \eqref{equ:upper_bound}, for any $x_0>0$, there exists $k_0>0$ such that 
$x_k\leq r+\varepsilon$ for all $k\geq k_0$. 
Since $f(x)$ is decreasing and $b(x)x$ is increasing with respect to $x$, we have 
\[ y_{k+1} \leq \big(s+b(r+\varepsilon)(r+\varepsilon)f(0)\big)y_k = \Big( s+\frac{cb_0(r+\varepsilon)}{1+\gamma(r+\varepsilon)} \Big)y_k, \]
for all $k\geq k_0$, 
which implies that $y_k\to 0$ as $k\to+\infty$. Thus, for any $\varepsilon>0$ satisfying $1-c\varepsilon >0$, there exists $k_1>k_0$ such that $y_k\leq \varepsilon$ for all $k\geq k_1$. Then one has
\[ x_{k+1} \geq \phi(x_k)x_k(1-f(0)y_k)\geq (1-c\varepsilon)\phi(x_k)x_k\ = e^{\ln (1-c\varepsilon)+r-x_k}x_k \]
for all $k\geq k_1$. 
Let $\{z_k\}$ be the sequence of the one dimensional Ricker model
\begin{equation*}
    z_{k+1}=e^{\ln (1-c\varepsilon)+r-z_k}z_k
\end{equation*}
with initial condition $z_{k_1}=x_{k_1}$. Since $\ln(1-c\varepsilon)+r\leq 1$, one has that $e^{\ln (1-c\varepsilon)+r-x}x$ is an increasing function from $[0,K_1]$ into $[0,K_1]$. Then 
by the similar arguments as the proof of Lemma \ref{lemma:leqs_Ricker} one can prove that
$
x_k\geq z_k
$ for all $k\geq k_1$. It follows that 
\begin{equation*}
    \liminf_{k\to +\infty} x_k\geq \lim_{k\to +\infty}z_k = r +\ln(1-c\varepsilon).
\end{equation*}
By letting $\varepsilon \to 0$ one has that
\begin{equation}\label{equ:lower_bound}
    \liminf_{k\to+\infty}x_k\ge r.
\end{equation}
As a consequence, $x_k\to r$ as $k\to +\infty$ by \eqref{equ:upper_bound} and \eqref{equ:lower_bound}. This completes the proof. 
\end{proof}

\subsection{The positive fixed point}\label{sec:positive-exist-local}

\begin{lemma}[Existence of the positive fixed point]\label{thm:positive-existence}
The system \eqref{pp_model} has a unique positive fixed point $p^*=(x^*,y^*)$ if and only if $r>0$ and $s+\frac{cb_0r}{1+\gamma r}>1$.
\end{lemma}
\begin{proof}
Assume that $p^*=(x^*,y^*)$ is a positive fixed point  of system \eqref{pp_model}. Then $(x^*,y^*)$ is a positive solution of the equations 
\begin{equation}\label{eq:positive-relation}
\begin{cases}
    \phi(x)(1-f(y)y)=1,\\
    s+b(x)x f(y)=1.
\end{cases}    
\end{equation}
Since $f(x)x<1$ for all $x\in [0,+\infty)$, one has 
$$
e^{r-x^*}=\phi(x^*)=\frac{1}{(1-f(y^*)y^*)}>1,
$$
which implies that $0<x^*<r$. By noticing that $b(x)x$ is an increasing function with respect to $x$, it then follows from \eqref{eq:positive-relation} that
$$
s+\frac{cb_0r}{1+\gamma r}=s+c b(r)r>s+c b(x^*)x^* >s+b(x^*)x^* f(y^*)=1.
$$
On the other hand, if $r>0$ and $s+\frac{cb_0r}{1+\gamma r}>1$, then there exists $\hat{x}\in(0,r)$ satisfying  $s+\frac{cb_0\hat{x}}{1+\gamma \hat{x}}=1$. 
 Define the functions $S(x)$ and $V(x)$ for $x\in[0,+\infty)$ as  
\begin{equation}\label{equ:S}
    S(x) = \frac{1}{c}\left( e^{r-x} -1\right)
\end{equation}
and 
\begin{equation}\label{equ:V}
    V(x) = \frac{1}{c}\big( \frac{cb_0}{1-s}\frac{x}{1+\gamma x}-1 \big).
\end{equation}
Note that the equations \eqref{eq:positive-relation} have a unique positive solution $p^*$ if and only if there exists $x^*\in (\hat{x},r)$ such that $S(x^*)=U(x^*)$. Clearly, 
$$
S(\hat{x})=\frac{1}{c}\left( e^{r-\hat{x}} -1\right) >0=U(\hat{x})
$$
and 
$$
S(r)=0<\frac{1}{c}\left( \frac{cb_0}{1-s}\frac{r}{1+\gamma r}-1 \right)=U(r).
$$
Therefore, there is a unique intersection of $S$ and $U$ because $S(x)$ is decreasing while $U(x)$ is increasing with respect to $x$. This completes the proof.
\end{proof}
\begin{remark}
Though there is a positive fixed point $p^*=(x^*, y^*)$ when $r>0$ and $s+\frac{cb_0r}{1+\gamma r}>1$, which are equivalent to 
\begin{equation}\label{equ:existence-pf}
    cb_0>(1-s)\gamma, \quad r>\frac{1-s}{cb_0-(1-s)\gamma},
\end{equation}
it is difficult to solve the explicit expressions of $x^*$ and $y^*$ with respect to the parameters of system \eqref{pp_model}. But it is not difficult to check that
\begin{equation}\label{eq:x-to-r}
    x^*=\zeta(r),
\end{equation}
where $\zeta(r)$ is an increasing function with the inverse $\zeta^{-1}$ given by
\begin{equation}\label{eq:r-to-x}
     \zeta^{-1}(x^*)= \ln{\Big(\frac{cb_0x^*}{(1-s)(1+\gamma x^*)}\Big)} +x^*
\end{equation}
and 
\begin{equation}\label{eq:y-to-r}
   y^*=\psi(r) :=\frac{b_0\zeta(r)}{(1-s)(1+\gamma \zeta(r))}-\frac{1}{c},
\end{equation}
that is,
\begin{equation}
    p^*= (\zeta(r),\psi(r)).
\end{equation}
\end{remark}

\begin{lemma}\label{lem:local-stability}
Assume that system \eqref{pp_model} has a positive fixed point $p^*$. Then $p^*$ is locally asymptotically stable if $r\leq 2+\ln\frac{2cb_0}{(1-s)(1+2\gamma)}$.
\end{lemma}
\begin{proof}
Let $J_{p^*}=DG(p^*)$ be the Jacobian matrix of $G$ at $p^*$ given by
\begin{equation}\label{eq:Jacbian-positive}
    J_{p^*}= \left(\begin{array}{cc}1-x^* & -\frac{cx^*}{1+cy^*}\\
    \frac{(1-s)y^*}{x^*(1+\gamma x^*)} & s+\frac{1-s}{1+cy^*} \end{array} \right)
\end{equation}
According to Jury condition (Lemma 4.14 in \cite{Niu2016ON} or \cite{allen2006introduction}), the positive fixed point $(x^*,y^*)$ is asymptotically stable if
\begin{equation}\label{equ:jury}
    |\mathrm{tr}~ J_{p^*}| < 1+\det J_{p^*} <2, 
\end{equation}
which is equivalent to 
\begin{subequations}\label{eq:Jury-abc}
\begin{equation}\label{eq:Jury-a}
    1+\det J_{p^*}+\mathrm{tr}~ J_{p^*} >0,
\end{equation}
\begin{equation}\label{eq:Jury-b}
    1+\det J_{p^*}-\mathrm{tr}~ J_{p^*} >0,
\end{equation}
\begin{equation}\label{eq:Jury-c}
    1-\det J_{p^*} >0,
\end{equation}
\end{subequations}
where $\mathrm{tr}~ J_{p^\ast}$ is the trace of $J_{p^\ast}$. 
By calculation, one has that
\[ 1+\det J_{p^*}-\mathrm{tr}~ J_{p^*} =\frac{(1-s)cy^\ast}{1+cy^\ast}\left(x^\ast+\frac{1}{1+\gamma x^\ast}\right)>0\]
and
\[ 1-\det J_{p^*} = sx^*+(1-s)\left( x^*+\frac{c\gamma x^* y^*}{1+\gamma x^*}\right)>0\]
because $0<s<1,x^*>0,y^*>0$, 
which implies that \eqref{eq:Jury-b} and \eqref{eq:Jury-c} hold. 
On the other hand, one has that
\begin{equation}\label{equ:jury-a-geq}
    \begin{aligned}
   1+\det J_{p^*}+\mathrm{tr}~ J_{p^*} &= (2-x^*)\big(1+s+\frac{1-s}{1+cy^*}\big)  + \frac{(1-s)cy^*}{(1+\gamma x^*)(1+cy^*)} \\
    &> (2-x^*)\big(1+s+\frac{1-s}{1+cy^*}\big).
\end{aligned}
\end{equation}
Therefore, if
\[ 
r \leq 2+\ln\frac{2cb_0}{(1-s)(1+2\gamma)},
\]
then one has $0<x^*\leq 2$ by \eqref{eq:x-to-r} and \eqref{eq:r-to-x},   
which ensures that \eqref{eq:Jury-a} holds by \eqref{equ:jury-a-geq}. 
\end{proof}

\begin{corollary}\label{cor:no-NS}
    Neimark-Sacker bifurcation does not occur for system \eqref{pp_model}.
\end{corollary}
\begin{proof}
    Note that the eigenvalues of $DG(0,0)$ and $DG(r,0)$ are real numbers, so the  Neimark-Sacker bifurcation cannot occur at the boundary fixed points. Suppose that system \eqref{pp_model} has a positive fixed point $p^*$. If the Jacobian $J_{p^*}=DG(p^*)$ has two complex roots, then  
    $
    0<\det J_{p^*}<1$, 
    since \eqref{eq:Jury-c} always hold. Therefore, Neimark-Sacker bifurcation cannot occur at a positive fixed point. This completes the proof.
\end{proof}

\subsection{Global stability of the positive fixed point}\label{sec:positive-global}
In the section, we study the global stability of the positive fixed point for system \eqref{pp_model} by using a geometric method. 
\begin{lemma}\label{lem:global-R}
 Assume that $S_i,i=1,2,$ are two strictly decreasing functions and $V$ is a strictly increasing function on $[0,+\infty)$, such that there exists $\omega \in (0,+\infty)$ satisfying 
$S_1(\omega) = S_2(\omega) = V(\omega)$. 
Suppose that $S_1(x) \geq S_2(x)$ for $x\in [0,\omega]$ and $S_1(x) \leq S_2(x) $ for $x\in [\omega,+\infty)$. 
Let $\mathcal{R}(S_i,\cdot)$ be the functions defined on some interval $[\omega_*,\omega^*]$ given by
\begin{equation}
    \mathcal{R}(S_i,y)=V\circ S_i^{-1} \circ V \circ S_i^{-1} (y),
\end{equation}
where $\omega^* =V(\omega)$ and $\omega_*\leq\omega^*$.
Then $\mathcal{R}(S_1,y)\geq \mathcal{R}(S_2,y)$ for all $y\in[\omega_*,\omega^*]$.
\end{lemma}
\begin{proof}
For any $y\in[\omega_*,\omega^*]$, one has that 
\begin{equation}\label{eq:R-ineq1}
   S_2^{-1}(y)\geq S_1^{-1}(y) \geq \omega
\end{equation}
by the properties of $S_1$ and $S_2$. It then follows from the monotonicity of $V$ and $ S_2^{-1}$ that
\begin{equation}\label{equ:VS}
    V\circ S_2^{-1}(y) \geq V \circ S_1^{-1}(y)\geq  \omega^*
\end{equation}
and furthermore,
\begin{equation}\label{eq:R-ineq3_1}
S_2^{-1}\circ V\circ S_1^{-1}(y) \geq S_2^{-1}\circ V\circ S_2^{-1}(y).
\end{equation}
Thus, by \eqref{equ:VS} and \eqref{eq:R-ineq3_1}, one has that
\begin{equation}\label{eq:R-ineq3_2}
   S_1^{-1}\circ V\circ S_1^{-1}(y)\geq  S_2^{-1}\circ V\circ S_1^{-1}(y)\geq S_2^{-1}\circ V\circ S_2^{-1}(y).
\end{equation}
Therefore, 
\begin{equation}\label{eq:R-ineq4}
    V\circ S_1^{-1}\circ V \circ S_1^{-1}(y)\geq V\circ S_2^{-1}\circ V\circ S_2^{-1}(y),
\end{equation}
that is,
\begin{equation}\label{eqineqR}
    \mathcal{R}(S_1,y)\geq \mathcal{R}(S_2,y).
\end{equation}
\end{proof}

\begin{theorem}[Global stability]\label{thm:global-stability}
Assume that system \eqref{pp_model} has a positive fixed point $p^*=(x^*,y^*)$, that is, condition \eqref{equ:existence-pf} holds. If $r\leq 1$ and 
\begin{equation}\label{equ:global-ineq}
    \gamma \zeta^2(r) + 2\zeta(r) > \frac{1-s+cb_0 r}{cb_0-(1-s)\gamma},
\end{equation}
then $p^*$ is globally asymptotically stable in $\dot{\mathbb{R}}^2_+$.
\end{theorem}
\begin{proof}
We use a modification of the method of Ackleh et al. \cite[Theorem 2.1]{ackleh2021nullcline} to prove the theorem. Recall the functions $S$ and $V$ give by \eqref{equ:S} and \eqref{equ:V}, respectively. Let $U$ be the function
defined on  $[0,+\infty)$ by
\begin{equation*}
    U(x) = \begin{cases}
    \begin{aligned}
        &- \frac{1}{r-x^*}\big( \frac{b_0 x^*}{(1-s)(1+\gamma x^*)}-\frac{1}{c} \big) (x-r), \quad &&x\in[0,r],\\
        &S(x), \quad &&x\in(r,+\infty).
    \end{aligned}
    \end{cases}
\end{equation*}
Obviously, $S$ and $U$ are strictly decreasing functions and $V$ is a strictly increasing function. Note that $$S(x^*)=U(x^*)=V(x^*)=y^*$$
and  
$$S(r)=U(r)=V(\hat{x})=0$$
 with $\hat{x} =\frac{1-s}{cb_0-\gamma(1-s)}$. Note that $S(x) > U(x)$ for $x\in [0, x^*)$ while $S(x) > U(x)$ for $x \in( x^*,+\infty)$.
Consider the functions $\mathcal{R}(S,\cdot)$ and $\mathcal{R}(U,\cdot)$ defined on $[y_*,y^*]$ given by
$$\mathcal{R}(S,y) =V\circ S^{-1} \circ V \circ S^{-1} (y)$$
and
$$\mathcal{R}(U,y) =V\circ U^{-1} \circ V \circ U^{-1} (y),$$
respectively, where
\begin{equation}\label{equ:y*}
    y_*= 
    \begin{cases}
    \begin{aligned}
    &0, \quad  &&\text{if~~}  V(r) < U (\hat{x}),\\
    &U \circ V^{-1} \circ U (\hat{x}), \quad &&\text{if~~} V(r) \geq U (\hat{x}).
    \end{aligned}
    \end{cases}
\end{equation}
It follows from Lemma \ref{lem:global-R} that
$$
\mathcal{R}(S,y) \geq \mathcal{R}(U,y) 
$$
for all $y\in[y_*,y^*]$. 
Note that if there exists $\tilde{y}\in(y_*,y^*)$ such that $\mathcal{R}(U,\tilde{y})=\tilde{y}$, then one has $\mathcal{R}(U,\bar{y})=\bar{y}$, where 
$\bar{y} = V\circ U^{-1}(\tilde{y})$ and $\bar{y}\in(y^*,U(\hat{x})]$. 
It is easy to check that the existence of solution of $\mathcal{R}(U,y)=y$ is equivalent to a quadratic function in $y$, which implies that there are at most two different solutions of $\mathcal{R}(U,y)=y$. 
It is not difficult to check that 
$$
\frac{d \mathcal{R}(U,y^*)}{dt} = \Big(\frac{cb_0(r-x^*)}{(1+\gamma x^*)( cb_0 x^* - (1-s)(1+\gamma x^* ) )} \Big)^2<1
$$
by \eqref{eq:x-to-r} and \eqref{equ:global-ineq}, so there exists $y_0\in (y_*,y^*)$ such that $\mathcal{R}(U,y_0)>y_0$. It follows that $\mathcal{R}(U,y)>y$ for all $y\in(y_*,y^*)$ by noticing that $\mathcal{R}(U,y^*)=y^*$. Thus, if $y_*=U \circ V^{-1} \circ U (\hat{x})$, then one has $y_*<\mathcal{R}(U,y_*)=0$, and hence $y_*\leq 0$ by \eqref{equ:y*}. Therefore, $\mathcal{R}(U,\cdot)$ can be defined on $[0,y^*]$ and $\mathcal{R}(U,y)>y$ for all $y\in[0,y^*)$, and hence, $\mathcal{R}(S,\cdot)$ can be defined on $[0,y^*]$ and 
\begin{equation}\label{equ:R(S,y)gey}
    \mathcal{R}(S,y)\geq \mathcal{R}(U,y) > y
\end{equation}
for all $y\in [0,y^*)$. 
Thus, $y^*$ is the unique fixed point of $\mathcal{R}(S,\cdot)$ on $[0,y^*]$. Let $A^m_0 = S^{-1}\circ V \circ S^{-1}(0)$, $A^M_0 = S^{-1}(0)= r$, $B^m_0=V(A^m_0)$ and $B^M_0 =V(A^M_0)$. It is clear that $x^*<A^M_0\leq K_1$ and $B^M_0>y^*$. By noticing that $\mathcal{R}(S,y)$ is strictly increasing and $ \mathcal{R}(S,0)> 0$, one has that $\hat{x}<A^m_0<x^*$ and $0<B^m_0<y^*$. Define the sets 
\[ \mathcal{D}_k = \{ (x,y)\in \mathcal{K}: A^m_k\leq x\leq A^M_k,~B^m_k\leq y\leq B^M_k \},~~k=1,2,\ldots\]
with 
\begin{equation}\label{equ:alpha}
    A^m_k = S^{-1}\circ V(A^M_{k-1}), ~~A^M_k = S^{-1}\circ V(A^m_{k-1})
\end{equation}
and
\begin{equation}\label{equ:beta}
    B^m_k = V\circ S^{-1}(B^M_{k-1}),~~B^M_k =V\circ S^{-1}(B^m_{k-1}).
\end{equation}
Note that $G_1(x,y)$ is increasing respect to $x\in[0,K_1]$ because $r\leq 1$ and $G_2(x,y)$ is increasing respect to $y$. By repeating the similar arguments as the proof of Theorem 2.1 in \cite{ackleh2021nullcline}, one can prove that each set $\mathcal{D}_k$ is positively invariant and absorbs all orbits of system $\eqref{pp_model}$ that originate in $(0,+\infty)^2$ by Lemma \ref{lem:absorption}. Since $S^{-1}\circ V$ and $V\circ S^{-1}$ are strictly decreasing, the sequences $\{A^m_k\},\{B^m_k\}$ are increasing and $\{A^M_k\},\{B^M_k\}$ are decreasing by \eqref{equ:alpha} and \eqref{equ:beta}, and moreover,
$$
0<A^m_0\leq A^m_k\leq x^*\leq A^M_k\leq r,~~0<B^m_0\leq B^m_k\leq y^*\leq B^M_k\leq V(r)
$$ 
for all $k\geq 1$. Therefore, $\mathcal{D}_{k+1}\subset \mathcal{D}_{k}$ for all $k\geq 1$ and $\{A^m_k\},\{B^m_k\},\{A^M_k\}$ and $\{B^M_k\}$ are convergent, and we denote their limits, respectively, by $A_*$, $B_*$, $A^*$ and $B^*$. Let 
$$\mathcal{D}=\{(x,y)\in \mathcal{K}: A_*\leq x\leq A^*,B_*\leq y \leq B^*\}.$$
Then by \eqref{equ:alpha} and \eqref{equ:beta}, one has that $S(A^*)=B_*$, $V(A_*)=B_*$, $V(A^*)=B^*$ and $\mathcal{R}(S,B_*)=B_*$, which implies that $B_*=B^*=y^*$ and $A_*=A^*=x^*$. Therefore, $D=\{ p^*\}$, and hence $p^*$ is globally asymptotically stable in $(0,+\infty)^2$. This completes the proof.
\end{proof}

\begin{corollary}\label{coro:global-stability}
Assume that the condition \eqref{equ:existence-pf} holds. Then the system \eqref{pp_model} has a globally asymptotically stable positive fixed point if  
\begin{equation}\label{equ:r0}
   \ln{ \Big( 1- \sqrt{\frac{cb_0-(1-s)\gamma}{cb_0(1+\gamma r )}}\Big)\frac{cb_0}{\gamma(1-s)}} + \frac{1}{\gamma}\big( \sqrt{\gamma \frac{1-s+cb_0 r}{c b_0 -(1-s)\gamma} +1}-1 \big)<r\leq 1.
\end{equation}
\end{corollary}
\begin{proof}
  Since the condition   \eqref{equ:global-ineq} holds if \eqref{equ:r0} is true, the conclusion follows from Theorem \ref{thm:global-stability} immediately. 
\end{proof}
\begin{example}\label{Exam:global}
Consider the system \eqref{pp_model} with the parameters $b_0=4,\gamma = \frac{3}{2},c=\frac{9}{10},s=\frac{1}{10}$, which is given by
\begin{equation}\label{equ:example}
    \begin{cases}    x_{k+1}=e^{r-x_k}\left( 1-\frac{9y_k}{10+9y_k}\right)x_k,\\
    y_{k+1}=\frac{y_k}{10}+\frac{72 x_k y_k}{(2+ 3 x_k)(10+9y_k)}.
\end{cases} 
\end{equation}
Note that $cb_0>(1-s)\left(\gamma+\frac{1}{2}\right)$, so the system \eqref{equ:example} has a positive fixed point if $r>\frac{2}{5}$ by \eqref{equ:existence-pf}. It is easy to check that the condition \eqref{equ:global-ineq} holds if $r\in (0.8184,1]$, and hence the system \eqref{equ:example} has a globally asymptotically stable positive fixed point by Corollary \ref{coro:global-stability}. See Fig. \ref{fig:global-stability} for the numerical simulations for $r=1$.
\begin{figure}[ht!]
    \centering    \includegraphics[width=0.68\textwidth]{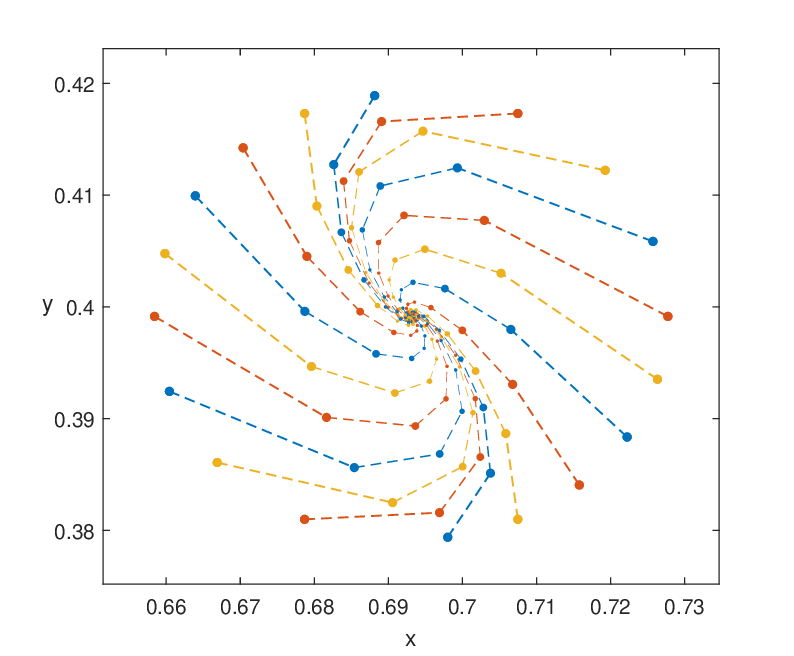}
    \caption{The orbits of the system \eqref{equ:example} through different initial values converge to the unique positive fixed point $(x^*,y^*)\approx (0.6930,0.3991)$. 
    }
    \label{fig:global-stability}
\end{figure}
    
\end{example}

\section{Period-doubling bifurcation}\label{sec:PD}
In this section, we study the occurrence of the period-doubling bifurcation in system \eqref{pp_model}. We first recall the period-doubling bifurcation theorem for the $n$-dimensional maps.

\begin{theorem}[Theorem 9.6 in \cite{Glendinning1994Stability}]\label{thm:demand-of-fild-stable}
Suppose that if $\mu=0$ the map $w_{k+1}=f\left(w_{k}, \mu\right), w \in \mathbb{R}^{n}$, has a fixed point at $w=0$ and that the linear map at $(w, \mu)=(0,0)$ has a simple eigenvalue of $-1$ and no other eigenvalues lie on the unit circle. Then the equation on the parameter-dependent center manifold is $x_{k+1}=H(x_{k}, \mu)$ with $H(0,0)=0$ and $\frac{\partial H}{\partial x}(0,0)=-1$. If
\[\sigma_1=2 H_{\mu x}+H_{\mu} H_{x x} \neq 0 ~~ \text{and}~~ \sigma_2=\frac{1}{2} H_{x x}^{2}+\frac{1}{3} H_{x x x} \neq 0, \]
(where the partial derivatives are evaluated at $(0,0)$), then a curve of periodic points of period two bifurcate from $(0,0)$ into $\mu>0$ if $\sigma_1 \sigma_2<0$ or $\mu<0$ if $\sigma_1 \sigma_2>0$. The fixed point from which these solutions bifurcate is stable in $\mu>0$ and unstable in $\mu<0$ if $\sigma_1>0$, with the signs of $\mu$ reversed if $\sigma_1<0$. The bifurcating cycle of period two is stable if it coexists with an unstable fixed point and vice versa. The bifurcation is said to be supercritical if the bifurcating solution of period two is stable and subcritical if the bifurcating solution is unstable.
\end{theorem}

\subsection{General outline of the calculations}\label{sec:General outline}

\begin{lemma}\label{lem:eigenvalue}
Assume that $cb_0>(1-s)\left(\gamma+\frac{1}{2}\right)$ and $r>\frac{1-s}{cb_0-(1-s)\gamma}$ so that system \eqref{pp_model} has a positive fixed point $p^*=(x^*,y^*)$. Then there exists $r^*>\frac{1-s}{cb_0-(1-s)\gamma}$ such that the two eigenvalues of $J_{p^\ast}=DG(p^*)$ are $-1$ and $-\det J_{p^\ast}$ with $-1<\det J_{p^\ast}<1$ when $r=r^*$.
\end{lemma}
\begin{proof}
Note that $-1$ is an eigenvalue of $J_{p^\ast}$ if and only if $$1+\mathrm{tr}~ J_{p^\ast}+\det J_{p^\ast}=0,$$
It follows from \eqref{eq:x-to-r} and \eqref{eq:y-to-r} and \eqref{eq:Jacbian-positive} that
\[ 
     1+\mathrm{tr}~ J_{p^\ast}+\det J_{p^\ast}
    =\Big(1+s+\frac{1-s}{1+c\psi(r)} \Big)(2-\zeta(r))+\frac{(1-s)c\psi(r)}{(1+\gamma \zeta(r))(1+c\psi(r))}.
\]
By the proof of Lemma \ref{lem:local-stability}, one has that
$$
    1+\mathrm{tr}~ J_{p^\ast}+\det J_{p^\ast}>0
$$
when $\frac{1-s}{cb_0-(1-s)\gamma}<r\leq 2+\ln\frac{2cb_0}{(1-s)(1+2\gamma)}$. 
Since 
$\zeta(r)\to +\infty$ and $\psi(r)\to \frac{b_0}{\gamma(1-s)}$ as $r\to +\infty$ by \eqref{eq:r-to-x} and \eqref{eq:y-to-r}, 
$$
    \underset{r\to +\infty}{\lim} (1+\mathrm{tr}~ J_{p^\ast}+\det J_{p^\ast}) = \underset{r\to +\infty}{\lim} \left(1+s+\frac{\gamma(1-s)^2}{cb_0}\right)(2-\zeta(r))=  -\infty.
$$
Therefore, there exists $r^*>\frac{1-s}{cb_0-(1-s)\gamma}$ such that $$1+\mathrm{tr}~ J_{p^\ast}+\det J_{p^\ast}=0$$
when $r=r^*$, and moreover, one eigenvalue of $J_{p^\ast}$ is $-1$ and the other eigenvalue is $-\det J_{p^\ast}$.
Together with \eqref{eq:Jury-b}, $1+\mathrm{tr}~ J_{p^\ast}+\det J_{p^\ast}=0 $ implies that 
$
1+\det J_{p^\ast}>0
$. 
Then together with \eqref{eq:Jury-c}, one has that 
$$
-1<\det J_{p^\ast}< 1
$$
when $r=r^*$. This completes the proof.
\end{proof}

Assume that $cb_0>(1-s)\left(\frac{1}{2}+\gamma\right)$ and $r\in (\frac{1-s}{cb_0-(1-s)\gamma},+\infty)$ such that the system \eqref{pp_model} has a positive fixed point $p^*=(x^*,y^*)=(\zeta(r),\psi(r))$. By Lemma \ref{lem:eigenvalue}, there exists $r^*>\frac{1-s}{cb_0-(1-s)\gamma}$ such that $J_{p^*}$ has an eigenvalue of $-1$ and its another eigenvalue $-\det J_{p^*}$ satisfying $-1<-\det J_{p^*}<1$  when $r=r^*$.  
We now include the
parameter $r$ as a new dependent variable as follows
\begin{equation}\label{Fild-firstsystem}
    \begin{cases}
        x_{k+1} = x_k e^{r_k-x_k} \big(1-\frac{cy_k}{1+cy_k}\big),\\
        r_{k+1} = r_k,\\
        y_{k+1} = s y_k + \frac{b_0 x_k}{1+\gamma x_k}\frac{c y_k}{1+ c y_k},
    \end{cases}
\end{equation}
with $r_k=r$, $k=1,2,\ldots$.
The system \eqref{Fild-firstsystem} has a positive fixed point $$\hat{p}^*_{r}=(x^*,r,y^*)=(\zeta(r),r,\psi(r)).$$
By making the change of variables given by $\tilde{x}_k=x_k-x^*$, $\tilde{r}_k=r_k-r^*$ and $\tilde{y}_k=y_k-y^*$, we transform system \eqref{Fild-firstsystem} into 
\begin{equation}\label{Fild-secondsystem}
    \begin{cases}
        \tilde{x}_{k+1}= \tilde{G}_1(\tilde{x}_k,\tilde{r}_k,\tilde{y}_k), \\
        \tilde{r}_{k+1}= \tilde{r}_k,\\
        \tilde{y}_{k+1}= \tilde{G}_2(\tilde{x}_k,\tilde{r}_k,\tilde{y}_k),
    \end{cases}
\end{equation}
so that the fixed point $\tilde{p}^*_{r^*}$ is translated to the origin $(0,0,0)$,
where
\[\begin{aligned}
    \tilde{G}_1(x,r,y)= &(x+\zeta(r+r^*)) e^{(r+r^*) - (x+\zeta(r+r^*))} \frac{1}{1+c(y+\psi(r+r^*))} -\zeta(r+r^*),\\
    \tilde{G}_2(x,r,y)= &s (y+\psi(r+r^*)) +\frac{b_0 (x+\zeta(r+r^*))}{1+\gamma (x+\zeta(r+r^*))}\frac{c(y+\psi(r+r^*))}{1+c(y+\psi(r+r^*))}- \psi(r+r^*).
\end{aligned}  \]
Note that the eigenvalues of the Jacobian matrix $\tilde{J}_{0}$ of system \eqref{Fild-secondsystem} at $(0,0,0)$ are $-1$, $1$ and $-\det J_{p^*}$. Thus, 
\[ T \tilde{J}_{0} T^{-1} = \left(\begin{array}{c}u_{k+1}\\ \mu_{k+1}\\v_{k+1} \end{array} \right) =  \left(\begin{array}{ccc} -1 & 0 & 0 \\ 0 & 1 & 0 \\ 0 & 0 &  -\det J_{p^*} \end{array}\right),  \]
where
\begin{equation*}
    T=\left(\begin{array}{ccc}
    \eta_1 & 0 & \eta_2 \\
    0 & 1 & 0 \\
    1 & 0 & 1 
    \end{array}\right),~~T^{-1}=\left(\begin{array}{ccc}
    -\eta_3 & 0 & \eta_4 \\
    0 & 1 & 0 \\
    \eta_3 & 0 & 1-\eta_4 
    \end{array}\right)
\end{equation*}
with
$$
\eta_1 = \frac{c\zeta(r^*)}{(2-\zeta(r^*))(1+c\psi(r^*))},~~
\eta_2=\frac{\zeta(r^*)(1+\gamma \zeta(r^*))(2-\zeta(r^*))}{(1-s)\psi(r^*)},$$
and $$
\eta_3 = \frac{1}{\eta_2-\eta_1}, ~~
\eta_4 = \frac{\eta_2}{\eta_2 - \eta_1}. 
$$
Let $\left(\begin{array}{c}\tilde{x}_k\\ \tilde{r}_k\\ \tilde{y}_k \end{array} \right)=T\left(\begin{array}{c}u_k\\ \mu_k\\v_k\end{array} \right)$. Then system \eqref{Fild-secondsystem} is equivalent to
\begin{equation}\label{Fild_Thirdsystem}
    \left(\begin{array}{c}u_{k+1}\\ \mu_{k+1}\\v_{k+1} \end{array} \right) =  \left(\begin{array}{ccc} -1 & 0 & 0 \\ 0 & 1 & 0 \\ 0 & 0 &  -\det J_{p^*} \end{array}\right) 
    \left(\begin{array}{c}u_k\\ \mu_k\\v_k 
    \end{array} \right) +\left(\begin{array}{c}L_1(u_k,\mu_k,v_k) \\ 0 \\L_3(u_k,\mu_k,v_k) 
    \end{array} \right),
\end{equation}
where the nonlinear terms
$$
L_1(u,\mu,v) =  u -\eta_3 \tilde{G}_1(\eta_1 u + \eta_2 v,\mu,u + v) + \eta_4 \tilde{G}_2(\eta_1 u + \eta_2 v,\mu,u + v)
$$
and
$$
L_3(u,\mu,v)  =  \det J_{p^*} v + \eta_3 \tilde{G}_1(\eta_1 u + \eta_2 v,\mu,u + v) +(1-\eta_4) \tilde{G}_2(\eta_1 u + \eta_2 v,\mu,u + v).
$$
Moreover, $L_1$ and $L_3$ can be expanded as
$$
    \begin{aligned}
        L_1(u,\mu,v)  =&~ \alpha_{200} u^2 + \alpha_{002}v^2 +\alpha_{110} u \mu + \alpha_{011} \mu v  +\alpha_{101} u v +\alpha_{300} u^3  +\alpha_{003} v^3\\
        & +\alpha_{210} u^2 \mu  + \alpha_{012} \mu v^2 + \alpha_{120} u \mu^2 + \alpha_{021} \mu^2 v +\alpha_{102} u v^2 + \alpha_{201} u^2 v  \\
        &+\alpha_{111} u \mu v + \text{h.o.t.}
        \end{aligned}
$$
and
$$
      \begin{aligned}   
        L_3(u,\mu,v)   =&~  \beta_{200} u^2 
        +\beta_{002}v^2 +\beta_{110} u \mu + \beta_{011} \mu v  + \beta_{101} u v +\beta_{300} u^3 + \beta_{003} v^3 \\
        & +\beta_{210} u^2 \mu  + \beta_{012} \mu v^2 + \beta_{120} u \mu^2 + \beta_{021} \mu^2 v +\beta_{102} u v^2 + \beta_{201} u^2 v  \\
        &+\beta_{111} u \mu v + \text{h.o.t.},\\ 
    \end{aligned}
$$
where h.o.t. denotes the higher order terms and 
\begin{equation}\label{equ:alpha-beta}
\begin{aligned}
    \alpha_{200}=&~\eta_4(\eta_1^2 j_{200} + \eta_1 j_{101} + j_{002}) -\eta_3(\eta_1^2 i_{200} + \eta_1 i_{101} + i_{002}) ,\\[2pt]
    \alpha_{110}=& ~\eta_4(\eta_1 j_{110} + j_{011})-\eta_3(\eta_1 i_{110} + i_{011}) ,\\[2pt]
    \alpha_{101}=&~\eta_4(2\eta_1\eta_2 j_{200} + \eta_1 j_{101} + \eta_2 j_{101} + 2j_{002}) \\
    &-\eta_3(2\eta_1\eta_2 i_{200} + \eta_1 i_{101} + \eta_2 i_{101} + 2i_{002}),\\[2pt]
    \alpha_{300}=&~\eta_4(\eta_1^3j_{300} + \eta_1^2 j_{201} + \eta_1 j_{102} + j_{003})\\
    &-\eta_3(\eta_1^3 i_{300} + \eta_1^2 i_{201} + \eta_1 i_{102} + i_{003})
    ,\\[2pt]
    \beta_{200}=&~\eta_3(\eta_1^2 i_{200} + \eta_1 i_{101} + i_{002}) + (1-\eta_4)(\eta_1^2 j_{200}+\eta_1 j_{101} + j_{002})
\end{aligned}
\end{equation}
with
\begin{equation}\label{equ:ij}
    i_{lmn}= \frac{1}{l!m!n!} \frac{\partial^{(l+m+n)} \tilde{G}_1(0,0,0)}{\partial{x}^l  \partial r^m  \partial{y}^n },\quad
    j_{lmn}= \frac{1}{l!m!n!} \frac{\partial^{(l+m+n)} \tilde{G}_2(0,0,0)}{\partial{x}^l  \partial{r}^m  \partial{y}^n }.
\end{equation}
From the center manifold theory \cite{Glendinning1994Stability,Marsden}, the system \eqref{Fild_Thirdsystem} can be reduced to a one-dimensional system by computing the parameter-dependent center manifold
\begin{equation}\label{eq:v_ApproximateExpansion}
\begin{aligned}
    v = h(u,\mu) = a_1 u^2+ a_2  u \mu +a_3 \mu^2+ a_4 u^3 + a_5 \mu^3 +a_6 u^2 \mu + a_7 u \mu^2 + \text{h.o.t.}. 
\end{aligned}
\end{equation}
Substituting \eqref{Fild_Thirdsystem} into \eqref{eq:v_ApproximateExpansion} gives 
$$
h(-u_k+L_1(u_k,\mu_k,h(u_k,\mu_k)),\mu_k)=-\det J_{p^*} h(u_k,\mu_k)+L_3(u_k,\mu_k,h(u_k,\mu_k)),
$$
from which one obtains that 
\begin{equation*}\label{eq:concrete_a}
    a_1= -\frac{\beta_{200}}{1+\det J_{p^*}}, \quad a_2= \frac{\beta_{110}}{1-\det J_{p^*}}, \quad
a_3= 0, \quad a_5=0.
\end{equation*}
Then the equation on the parameter-dependent center manifold is given by $u_{k+1}=H(u_k,\mu)$ with
\begin{equation*}\label{eq:Fild_H}
    H(u,\mu) = -u + L_1(u,\mu,h(u,\mu))\\
    = (-1 +d_1 \mu)u + 
    d_2u^2 + d_3 u^3+ \text{h.o.t.},
\end{equation*}
where 
\begin{equation*}
    d_1 = \alpha_{110}, \quad d_2 = \alpha_{200},\quad d_3 =\alpha_{300}- \frac{\alpha_{101}\beta_{200}}{1+\det J_{p^*}}.
\end{equation*}
Obviously, one has that $H(0,0)=0$ and $\frac{\partial H(0,0)}{\partial u} =-1$. 
Let 
\begin{equation}\label{eq:pd_fildstable1}
     \sigma_1 = 2 H_{\mu u}+H_{\mu} H_{u u}
        = 2 d_1 = 2 \alpha_{110}, 
\end{equation}
and 
\begin{equation}\label{eq:pd_fildstable2}
    \sigma_2 = \frac{1}{2} H_{u u}^{2}+\frac{1}{3} H_{u u u}
    =  \frac{1}{2} d_2^2 +\frac{1}{3} d_3 = \frac{1}{2} \alpha_{200}^2 +\frac{1}{3} \big( \frac{\alpha_{300}-\alpha_{101}\beta_{200}}{1+\det J_{p^*}} \big),
\end{equation}
where the partial derivatives of $H$ are evaluated at $(0,0)$.
Then by Theorem \ref{thm:demand-of-fild-stable}, we obtain the following theorem.
\begin{theorem}[Period-doubling bifurcation]\label{thm:fild-stable}
Assume that 
$cb_0>(1-s)\left(\gamma+\frac{1}{2}\right)$ and $r>\frac{1-s}{cb_0-(1-s)\gamma}$ so that system \eqref{pp_model} has a positive fixed point $p^*$. Then there exists $r^*>\frac{1-s}{cb_0-(1-s)\gamma}$ such that the Jacobian matrix $J_{p^*}$ of the system \eqref{pp_model} at the positive fixed point $p^*$ has an eigenvalue of $-1$, and the following conclusions hold.
\begin{enumerate}
    \item[(i)] If $\sigma_1<0$ and $\sigma_2>0$, then a supercritical period-doubling bifurcation occurs at $r=r^*$, and there exists a $\delta>0$ such that system \eqref{pp_model} has a stable 2-periodic orbit for $r\in (r^*,r^*+\delta)$, where $\sigma_1$ and $\sigma_2$ are defined by \eqref{eq:pd_fildstable1} and \eqref{eq:pd_fildstable2}.
    \item[(ii)] If $\sigma_1>0$ and $\sigma_2>0$, then a supercritical period-doubling bifurcation occurs at $r=r^*$, and there exists a $\delta>0$ such that system \eqref{pp_model} has a stable 2-periodic orbit for $r\in (r^*-\delta, r^*)$.
    \item[(iii)] If $\sigma_1>0$ and $\sigma_2<0$, then a subcritical period-doubling bifurcation occurs at $r=r^*$, and there exists a $\delta>0$ such that system \eqref{pp_model} has an unstable 2-periodic orbit for $r\in (r^*,r^*+\delta)$.
    \item[(iv)] If $\sigma_1<0$ and $\sigma_2<0$, then a subcritical period-doubling bifurcation occurs at $r=r^*$, and there exists a $\delta>0$ such that system \eqref{pp_model} has an unstable 2-periodic orbit for $r\in (r^*-\delta, r^*)$.
\end{enumerate}
\end{theorem}

\subsection{Period-doubling cascades to chaos}\label{sec:chaos}
In this subsection, we study the period-doubling bifurcations in system \eqref{pp_model} by analyzing the concrete system \eqref{equ:example} with a free parameter $r$ given in Example \ref{Exam:global}. 
Recall that if $r>\frac{2}{5}$, the system \eqref{equ:example} has a positive fixed point given by
$$p^*=(x^*,y^*)=(\zeta(r),\varphi(r))$$
with 
\[ \zeta^{-1}(x^*)= \ln{\big(\frac{8 x^*}{2+ 3 x^*}\big)} +x^*\]
and 
\[\psi(r)=\frac{80\zeta(r)}{9(2+3 \zeta(r))}-\frac{10}{9}. \]
By Lemma \ref{lem:eigenvalue}, for
\begin{equation}\label{equ:r*}
    r= r^*\approx 2.7732,
\end{equation}
the Jacobian matrix $J_{p^*}$ has an eigenvalue of $-1$ and the other eigenvalue $-\det J_{p^*}\approx 0.4746$. By including the parameter $r$ as a new dependent variable and translating the positive fixed point to the origin when $r=r^*$, the system \eqref{equ:example} is equivalent to
\begin{equation*}
    \begin{cases}    x_{k+1}=x_ke^{r_k-x_k}\big( 1-\frac{9y_k}{10+9y_k}\big),\\
    r_{k+1} = r_k,\\
    y_{k+1}=\frac{1}{10}y_k+\frac{72 x_k y_k}{(2+3 x_k)(10+9y_k)}. 
\end{cases} 
\end{equation*}
Let 
\[\begin{aligned}
    \tilde{G}_1(x,r,y)= &(x+\zeta(r+r^*)) e^{(r+r^*) - (x+\zeta(r+r^*))} \frac{10}{10+9(y+\psi(r+r^*))} -\zeta(r+r^*),\\
    \tilde{G}_2(x,r,y)= & \frac{1}{10} (y+\psi(r+r^*)) +\frac{72(x+\zeta(r+r^*))(y+\psi(r+r^*))}{ \big(2+3 (x+\zeta(r+r^*))\big) \big(10+9(y+\psi(r+r^*))\big)} - \psi(r+r^*).
\end{aligned}  \]
Then one has
$$
\sigma_1=2\alpha_{110}\approx -1.9363<0
$$
and
$$
\sigma_2=\frac{1}{2} \alpha_{200}^2 +\frac{1}{3} \big( \frac{\alpha_{300}-\alpha_{101}\beta_{200}}{1+\det J_{p^*}} \big)\approx 9.7182>0,
$$
where $\alpha_{110}$, $\alpha_{200}$, $\alpha_{101}$, $\alpha_{300}$ and $\beta_{200}$ are defined by \eqref{equ:alpha-beta}.
Therefore, by Theorem \ref{thm:fild-stable} (i) a supercritical period-doubling bifurcation occurs at $r=r^*$ given by \eqref{equ:r*}, and there exists a $\delta>0$ such that system \eqref{equ:example} has a stable 2-periodic orbit for $r\in (r^*,r^*+\delta)$. 

\begin{figure}[ht!]
    \centering
    \includegraphics[width =0.92\textwidth]{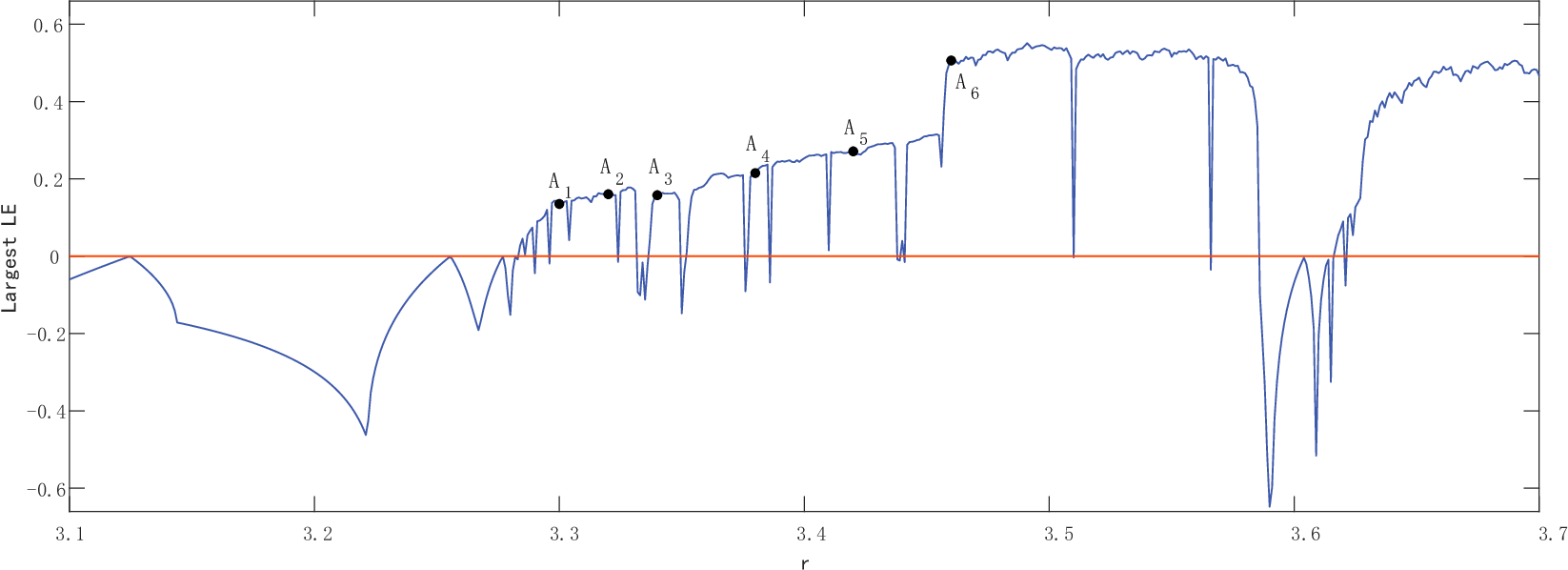}
    \caption{The largest Lyapunov exponent  as the function of $r$.
    The parameter values used in Fig. \ref{fig:chaos} with the corresponding largest Lyapunov exponents are marked with the closed dots $\bullet$,  where $A_1=(3.3,0.1352)$, $A_2=(3.32,0.1604)$, $A_3=(3.34,0.1578)$, $A_4=(3.38,0.2150)$, $A_5=(3.42,0.2712)$ and $A_6=(3.46,0.5063)$. 
    }
    \label{fig:Lyapunov exponent}
\end{figure}

\begin{figure}[ht!]
\centering
\begin{subfigure}{0.4\textwidth}
\includegraphics[width=\linewidth]{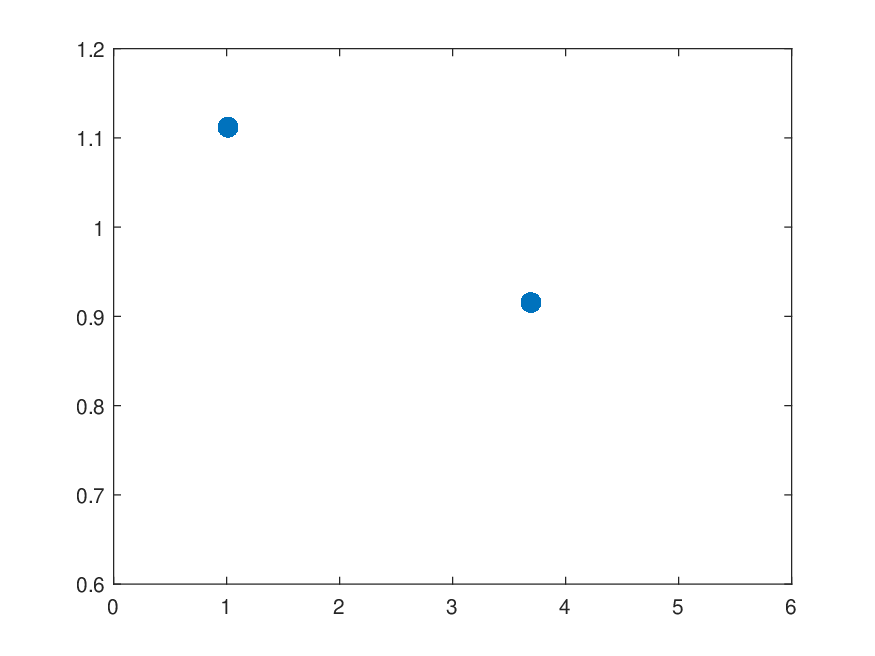}
\subcaption{$r=3.00$}
\label{fig:r=3.00}
\end{subfigure}
\begin{subfigure}{0.4\textwidth}
\includegraphics[width=\linewidth]{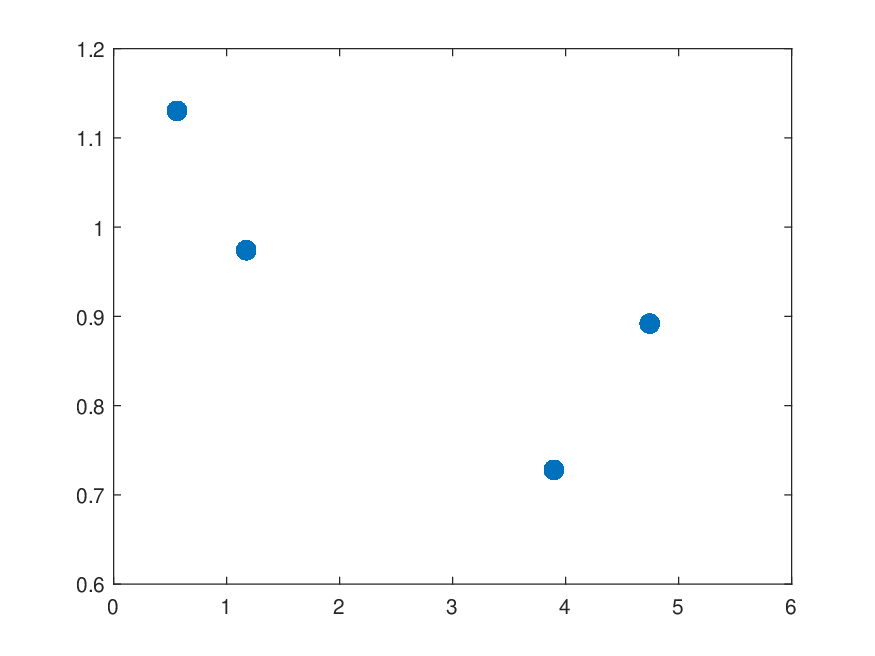}
\subcaption{$r=3.20$}
\label{fig:r=3.20}
\end{subfigure}\\
\begin{subfigure}{0.4\textwidth}
\includegraphics[width=\linewidth]{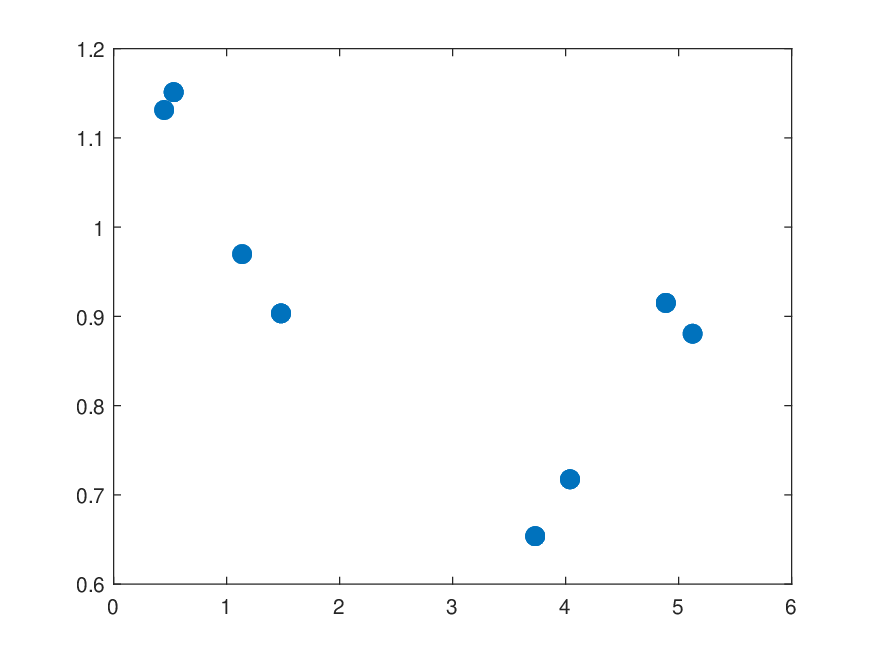}
\subcaption{$r=3.27$}
\label{fig:r=3.27}
\end{subfigure}
\begin{subfigure}{0.4\textwidth}
\includegraphics[width=\linewidth]{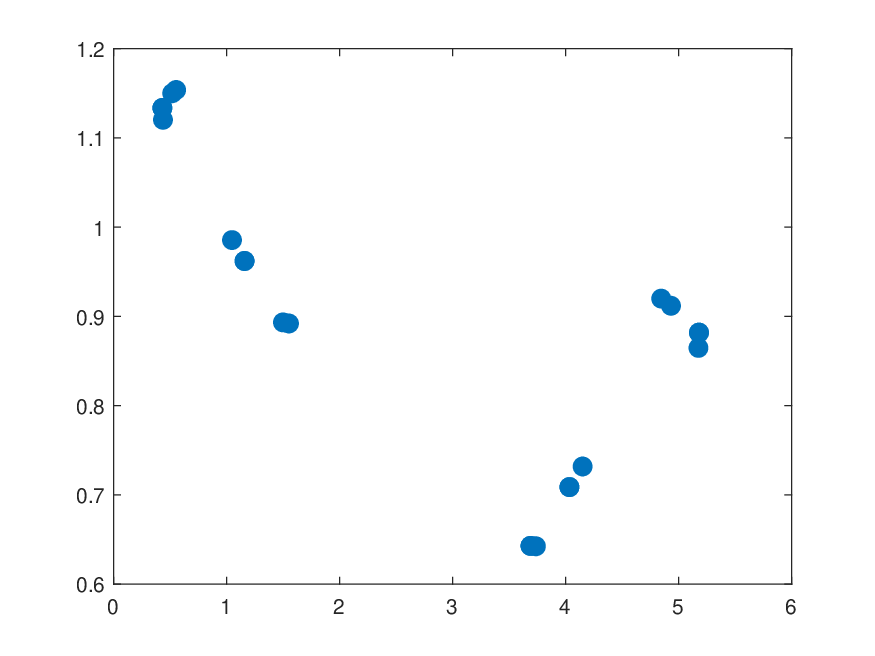}
\subcaption{$r=3.28$}
\label{fig:r=3.28}
\end{subfigure}
\caption{The different periodic attractors of the system \eqref{equ:example}. The initial value used is $(x_0,y_0)=(1,1)$. (a) The attractor at $r = 3.00$ is a period-2 cycle. (b) The attractor at $r = 3.20$ is a period-4 cycle. (c) The attractor at $r = 3.27$ is a period-8 cycle.  (d) The attractor at $r = 3.28$ is a period-16 cycle.}
\label{fig:period doubling}
\end{figure}

\begin{figure}[ht!]
\centering
\begin{subfigure}{0.4\textwidth}
\includegraphics[width=\linewidth]{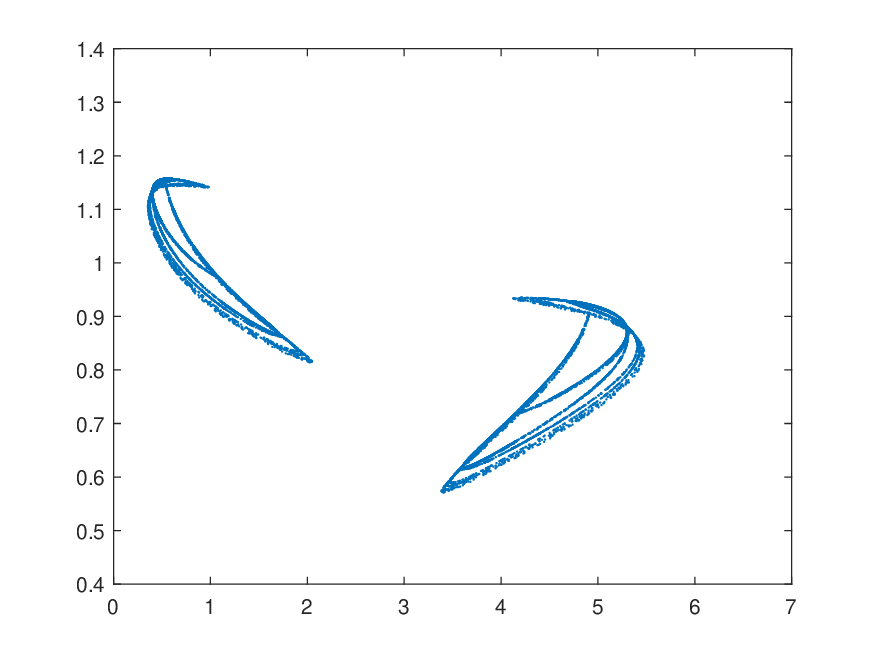}
\subcaption{$r=3.30$}
\label{fig:r=3.30}
\end{subfigure}
\begin{subfigure}{0.4\textwidth}
\includegraphics[width=\linewidth]{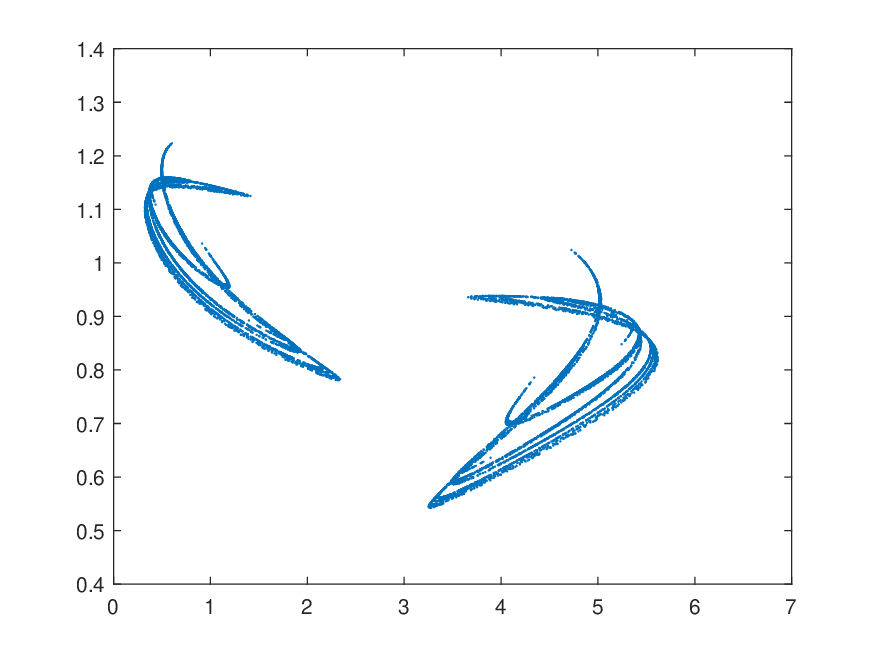}
\subcaption{$r=3.32$}
\label{fig:r=3.32}
\end{subfigure}\\
\begin{subfigure}{0.4\textwidth}
\includegraphics[width=\linewidth]{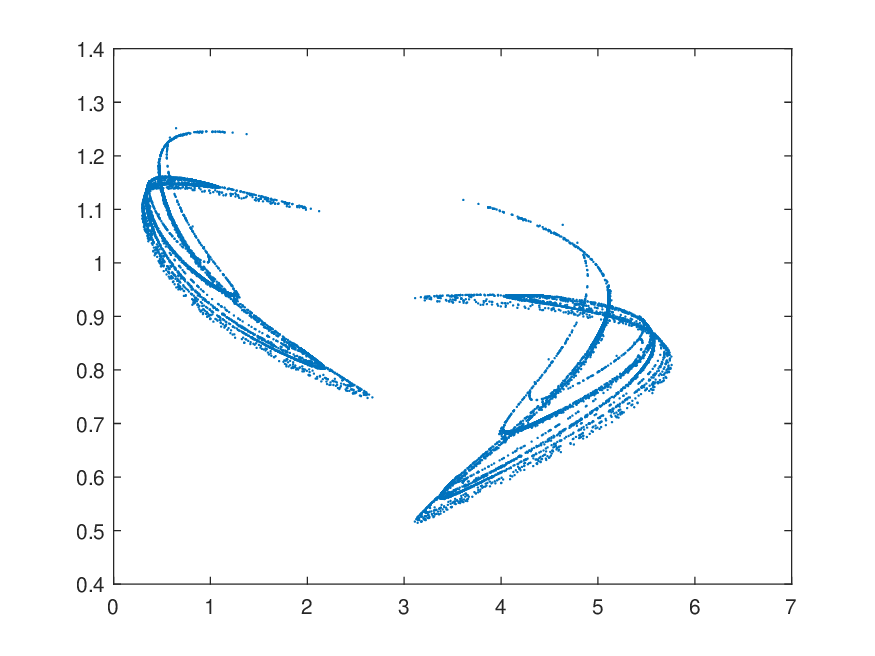}
\subcaption{$r=3.34$}
\label{fig:r=3.34}
\end{subfigure}
\begin{subfigure}{0.4\textwidth}
\includegraphics[width=\linewidth]{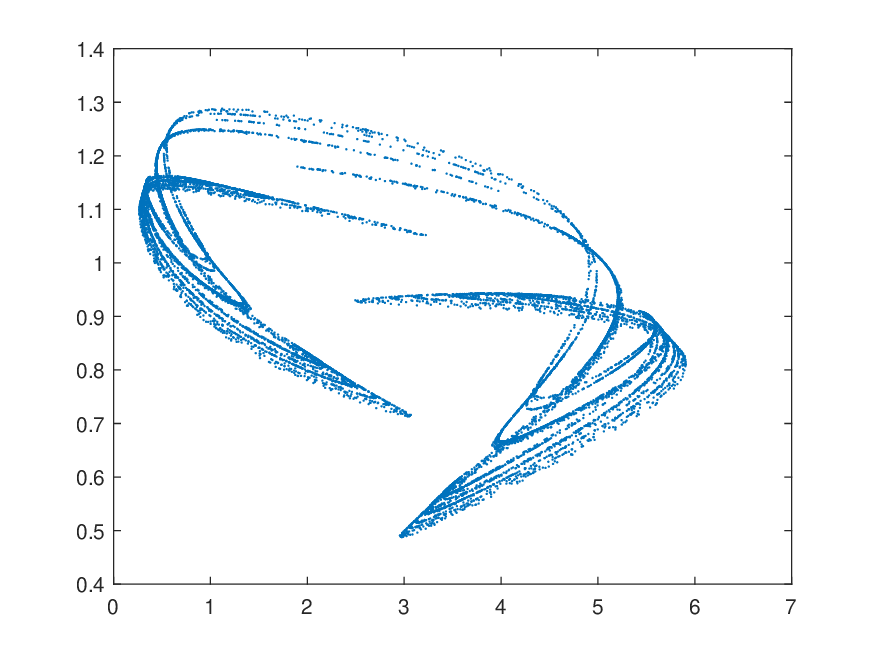}
\subcaption{$r=3.36$}
\label{fig:r=3.36}
\end{subfigure}\\
\begin{subfigure}{0.4\textwidth}
\includegraphics[width=\linewidth]{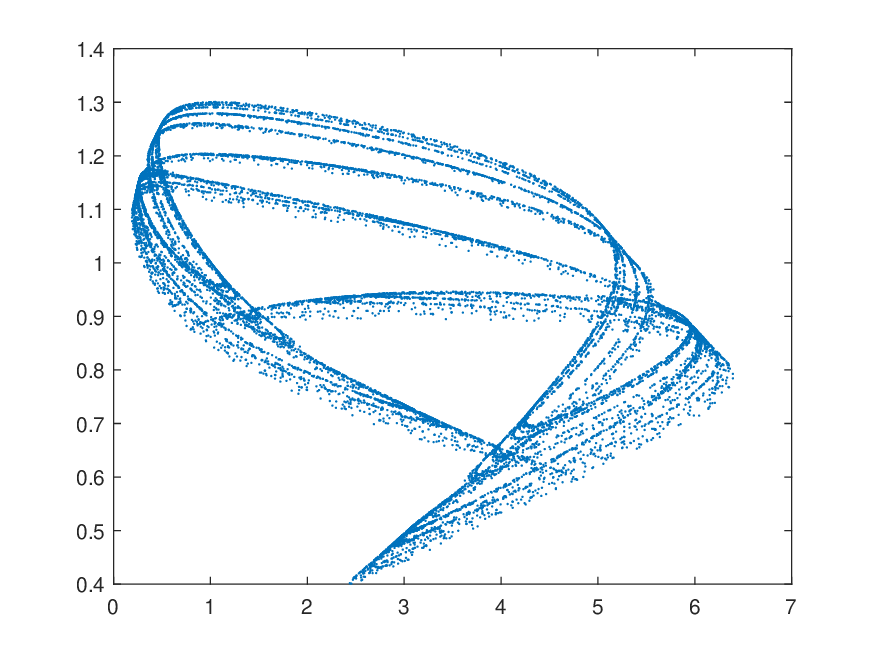}
\subcaption{$r=3.42$}
\label{fig:r=3.42}
\end{subfigure}
\begin{subfigure}{0.4\textwidth}
\includegraphics[width=\linewidth]{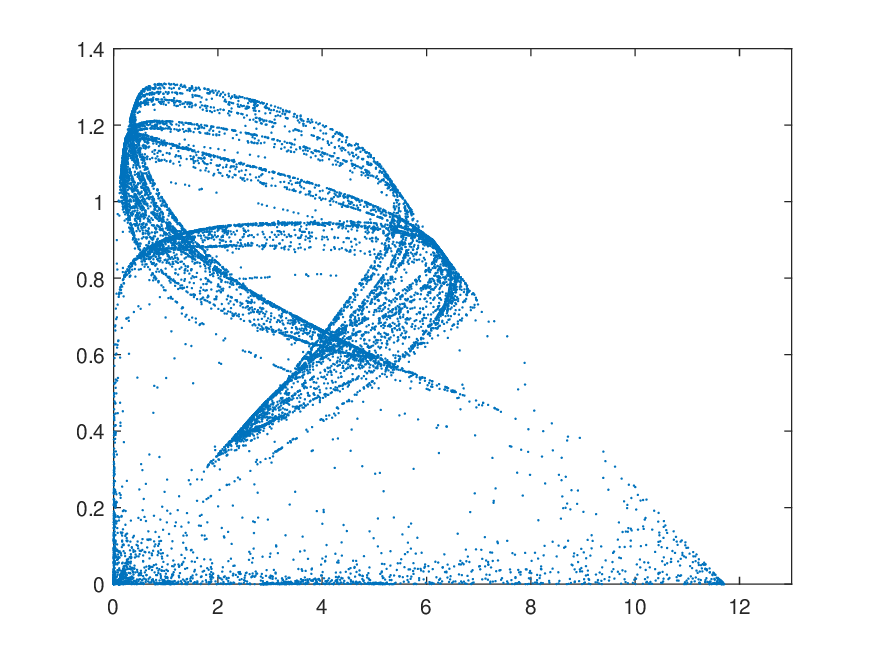}
\subcaption{$r=3.46$}
\label{fig:r=3.46}
\end{subfigure}
\caption{The different chaotic attractors of the system \eqref{equ:example}. The initial value used is $(x_0,y_0)=(1,1)$.}
\label{fig:chaos}
\end{figure}

\begin{figure}[ht!]
\centering
\begin{subfigure}{0.4\textwidth}
\includegraphics[width=\linewidth]{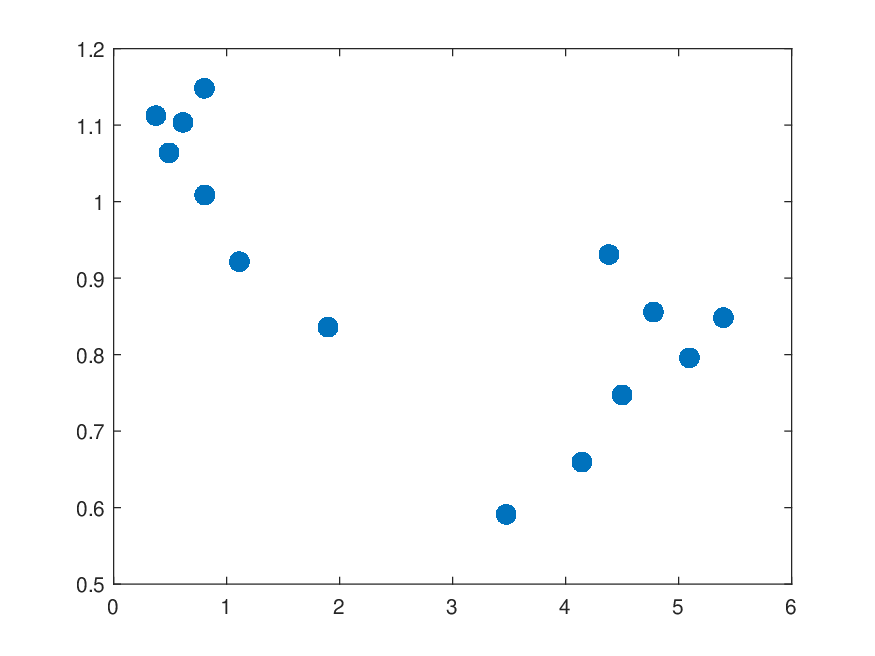}
\subcaption{$r=3.29564$}
\label{fig:14-periodic}
\end{subfigure}
\begin{subfigure}{0.4\textwidth}
\includegraphics[width=\linewidth]{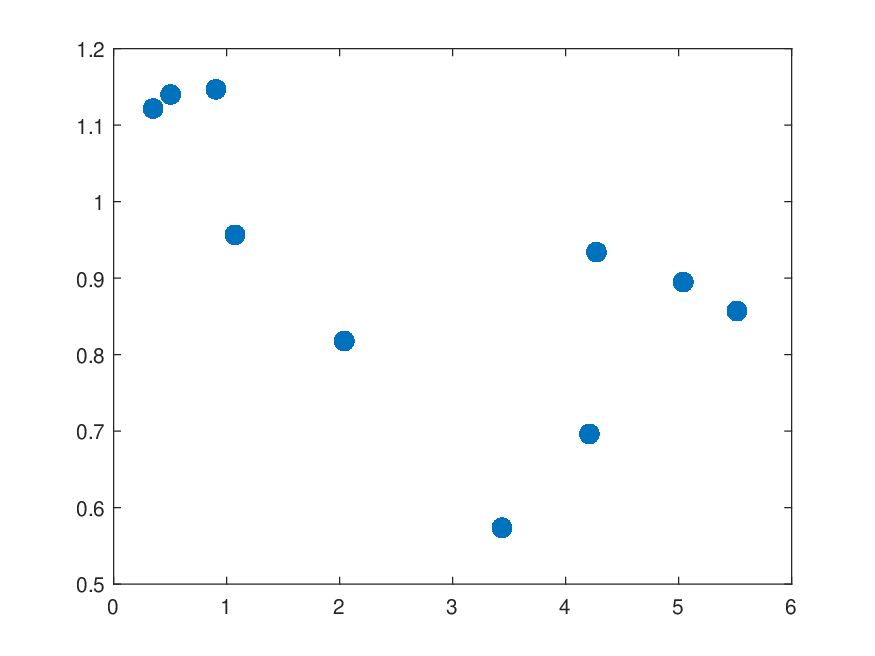}
\subcaption{$r=3.33142$}
\label{fig:10-periodic}
\end{subfigure}\\
\begin{subfigure}{0.4\textwidth}
\includegraphics[width=\linewidth]{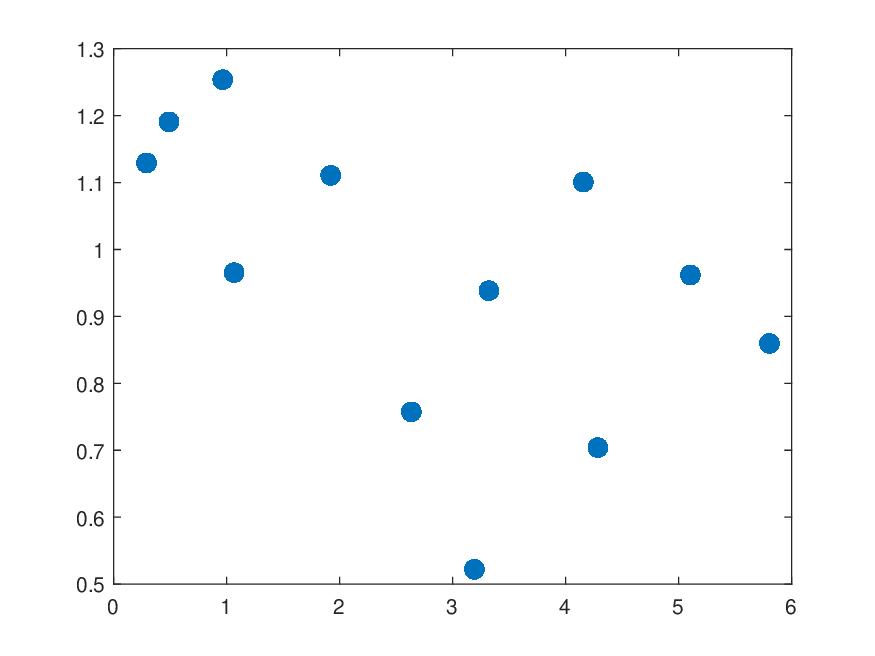}
\subcaption{$r=3.38593$}
\label{fig:12-periodic}
\end{subfigure}
\begin{subfigure}{0.4\textwidth}
\includegraphics[width=\linewidth]{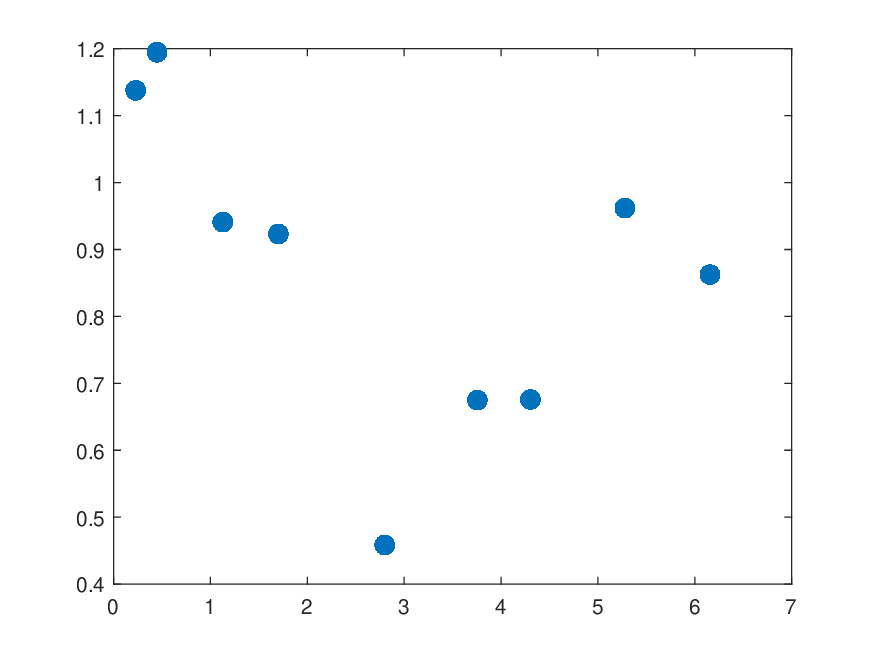}
\subcaption{$r=3.43853$}
\label{fig:9-periodic}
\end{subfigure}
\caption{The different periodic attractors of the system \eqref{equ:example}. The initial value used is $(x_0,y_0)=(1,1)$. (a) The attractor at
$r=3.29564$ is a period-14 cycle. (b) The attractor at
$r=3.33142$ is a period-10 cycle. (c) The attractor at
$r=3.38593$ is a period-12 cycle.  (d) The attractor at
$r=3.43853$ is a period-9 cycle. }
\label{fig:periodic_orbit}
\end{figure}

We now study the evolution of the dynamics of system \eqref{equ:example} as $r$ is increased from $r^*$ by using numerical simulations.  We calculate the largest Lyapunov exponent by using MatContM (\cite{Dhooge2008,Neirynck2016}) and find that it can be positive for some $r>3.30$ (see Fig. \ref{fig:Lyapunov exponent}), which indicates the occurrence of chaotic attractors (\cite{Alligood2000}). The route to chaos is of particular interest. As shown by our numerical simulations, a cascade of period-doubling bifurcations occur as $r$ is increased from $r^*$. First, the attracting 2-periodic orbit appears as  $r$ is increased from $r^*$ until about $r=3.1247$.  Subsequently, the attracting 4-period orbit appears until  $r>3.1247$. Then there are the attracting $8$-periodic orbit in sequence until $r>3.2555$,
and the attracting $16$-periodic orbit until $r>3.2770$, and so on.  Figs. \ref{fig:r=3.00}--\ref{fig:r=3.28} show the attracting 2-periodic orbit, 4-periodic orbit, 8-periodic orbit and 16-periodic orbit at $r=3$, $r=3.2$, $r=3.27$ and $r=3.28$ respectively. Until $r>3.2836$, the chaotic orbit occurs. See Fig. \ref{fig:chaos} for the different chaotic attractors for different values of the parameter $r$. Note that there are many periodic windows when $r>3.2836$ by Fig. \ref{fig:Lyapunov exponent}. For example, the 9-periodic orbit, 10-periodic orbit, 14-periodic orbit, 24-periodic orbit are detected as shown in Fig. \ref{fig:9-periodic}--\ref{fig:14-periodic}.

\section{Discussion}\label{sec:conclusion}
In this paper, we study the dynamics of the discrete Kolmogorov predator-prey model 
\begin{equation}\label{pp_model-1}
    \begin{cases}
        x_{k+1}=x_ke^{r-x_k}(1-\frac{cy_k}{1+cy_k}),\\
        y_{k+1}=sy_k+\frac{b_0x_k}{1+\gamma x_k}\frac{cy_k}{1+cy_k},
    \end{cases}
\end{equation}
of which the prey obeys the Ricker-type growth in the absence of predator, 
where $b_0,\gamma>0$, $0<c<1$ and $0<s<1$. By analyzing the stability of the positive fixed point, we find that the system exhibits period-doubling bifurcations at the positive fixed point. Using the center manifold theory, we are able to determine the conditions under which these bifurcations occur. We also find the route to chaos by using numerical simulations, that is,
period-doubling cascades lead to the chaos. This means that in addition to stable fixed points and periodic orbits, the system \eqref{pp_model-1} can possess chaotic attractors. Such chaotic behavior may stem from the interactions between species and the feedback loops within the ecosystem. This is a significant departure from the behavior of continuous predator-prey systems, which typically exhibit stable steady states or limit cycles. The discrete predator-prey system \eqref{pp_model-1} will provide
a reasonable explanation for those ecological communities in which populations are
observed to oscillate in a rather reproducible chaotic manner. Additionally, we have demonstrated the possibility of a globally asymptotically stable positive fixed point in system \eqref{pp_model-1} under certain conditions, indicating the potential for coexistence of prey and predator in a stable steady state. It should be pointed out that the complexity of the system presents challenges in providing a global analysis of the system, such as global stability of the fixed points or periodic orbits. Though we do not address the global stability of the periodic orbits for system \eqref{pp_model-1} in this paper, the geometric method provided here also offers a promising approach for future research on the global stability of periodic points.

\bibliographystyle{siamplain}

\end{document}